\documentclass[12pt]{article}

\usepackage{amsmath,amssymb,amsfonts,theorem,makeidx,latexsym,epsfig,,subfigure}

\newtheorem{defn}{Definition}[section]

\newtheorem{lemma}[defn]{Lemma}

{\theorembodyfont{\rmfamily}

\newtheorem{ex}[defn]{Example}}

\newtheorem{thm}[defn]{Theorem}

\newtheorem{prop}[defn]{Proposition}

\newtheorem{cor}[defn]{Corollary}

\numberwithin{equation}{section}

\newcommand{\h}{{\cal H}}

\newcommand{\ltr}{ L^2(\mathbb R) }

\newcommand{\ltn}{{\ell}^2(\mathbb N)}

\newcommand{\ltz}{{\ell}^2(\mathbb Z)}

\newcommand{\lti}{{\ell}^2(I)}

\newcommand{\si}{S^{-1}}

\newcommand{\mn}{\mathbb N}

\newcommand{\mr}{\mathbb R}

\newcommand{\mz}{\mathbb Z}

\newcommand{\mts}{ \{E_{mb}T_{na}g \}_{m,n \in \mz}}

\def\bp{{\noindent\bf Proof. \ }}

\def\ep{\hfill$\square$\par\bigskip}

\def\bqs{\begin{equation}}

\def\eqs{\tag*{$\square$}\end{equation}\par\bigskip}

\def\la{\langle}

\def\ra{\rangle}

\def\ga{\gamma}

\def\ftk{\{f_k\}_{k=1}^\infty}

\def\ctz{\{c_k\}_{k\in \mz}}

\def\etk{\{e_k\}_{k=1}^\infty}

\def\nl{\left|\left|}

\def\nr{\right|\right|}

\def\span{\overline{\text{span}}}

\def\Span{\text{span}}

\def\supp{\text{supp}}

\def\vn{\vspace{.1in}\noindent}

\def\bop{\begin{op}\rm}

\def\eop{\end{op}}

\def\cra{{\cal R}}

\def\bee{\begin{eqnarray}}

\def\ene{\end{eqnarray}}

\def\bes{\begin{eqnarray*}}

\def\ens{\end{eqnarray*}}

\def\bei{\begin{itemize}}

\def\eni{\end{itemize}}

\def\bt{\begin{thm}}

\def\et{\end{thm}}

\def\bc{\begin{cor}}

\def\ec{\end{cor}}

\def\bpr{\begin{prop}}

\def\epr{\end{prop}}

\def\bl{\begin{lemma}}

\def\el{\end{lemma}}

\def\bd{\begin{defn}}

\def\ed{\end{defn}}

\def\bex{\begin{ex}}

\def\enx{\end{ex}}

\def\bfi{\begin{fig}}

\def\efi{\end{fig}}

\def\inr{\int_{-\infty}^\infty}

\def\etz{\{e_k\}_{k\in \mz}}
\def\ftz{\{f_k\}_{k\in \mz}}
\def\ftki{\{f_k\}_{k\in I}}
\def\ctki{\{c_k\}_{k\in I}}
\def\sumk{\sum_{k\in I}}
\def\Tni{\{T^k \varphi\}_{k\in I}}
\def\ftz{\{f_k\}_{k\in \mz}}

\def\eptz{\{\varepsilon_k\}_{k\in \mz}}
\newcommand{\Tnz}{\{T^kf_0\}_{k\in\mz}}

\def\eptk{\{\varepsilon_k\}_{k=1}^\infty}

\title{
An open problem concerning operator representations of frames}

\date{\today}

\author{Ole Christensen, Marzieh Hasannasab }

\begin{document}

\maketitle

\begin{abstract} Recent research has shown that the properties of
overcomplete
Gabor frames and frames arising from shift-invariant systems form a precise match with certain conditions
that are necessary for a frame in $\ltr$ to have a representation $\{T^k \varphi\}_{k=0}^\infty$
for some bounded linear operator $T$ on $\ltr$ and some $\varphi \in \ltr.$ However,
for frames of this type  the existence of such
a representation has only been confirmed in the case of
Riesz bases. This leads
to several open questions connecting dynamical sampling, coherent states, frame
theory, and operator theory. The key questions can either be considered in the
general functional analytic context of operators on a Hilbert space, or
in the specific situation of Gabor frames in $\ltr.$
\end{abstract}

\section{Introduction and motivation}

A coherent state is a (typically overcomplete) system of vectors in a Hilbert space $\h.$
In general it is given by the action of a class of linear operators on  a single
element in the underlying Hilbert space. In particular, it could be given by
iterated action of a fixed operator on a single element, i.e., as
$\{T^k \varphi\}_{k=0}^\infty$ for some $\varphi \in \h$ and a
linear operator $T: \h \to \h.$ Coherent states play an important role in
mathematical physics \cite{Klauder,Perelomov}, operator theory, and modern harmonic analysis
\cite{Da1,CB}. In particular, a Gabor system (see the definition below) is a coherent
state.

Systems of vectors on the form $\{T^k \varphi\}_{k=0}^\infty$ also appear in
the more recent context of dynamical sampling \cite{A1,A2,A3,FP}. Here the goal is to analyse the
frame properties of such systems for certain classes of operators, e.g., normal
operators or self-adjoint operators. A different approach was taken in the
papers \cite{olemmaarzieh-1,olemmaarzieh-2}: here, the starting point is a frame
and the question is when and how it can be represented on the form
\bee \label{70404a} \{T^k \varphi\}_{k=0}^\infty \mbox{ for some } \varphi \in \h
 \mbox{ and a bounded linear
operator } T.\ene In particular the existence of such a representation was proved
for the case where the frame is indeed a basis, so the following discussion deals with
the overcomplete case.
It turns out that the conditions for the desired representation
to exist form a perfect match with the known
properties of Gabor frames in $\ltr.$ This raises the natural question whether
(some or all) overcomplete Gabor frames indeed have a representation
of the form \eqref{70404a}. A positive answer to the question would shed new light on Gabor frames, and also provide
new insight in the context of dynamical sampling. In fact, only few overcomplete
frames are known to have a representation on the form \eqref{70404a}. To the best knowledge of the authors, Gabor frames have not
been considered in this context before, except for the examples appearing
in the papers \cite{olemmaarzieh-1,olemmaarzieh-2}.

The above questions will also be analyzed with a different indexing, i.e., considering
systems on the form $\{T^k \varphi\}_{k=-\infty}^\infty$ instead of  $\{T^k \varphi\}_{k=0}^\infty.$ The indexing in terms of $\mz$ is natural for several well-known
classes of frames, and the theoretical conditions for a frame having such a representation
with a bounded operator $T$ are similar to the ones  for systems indexed by $\mn_0.$
However, the change in indexing gives an interesting twist on the problem: indeed, a
shift-invariant system always has a representation $\{T^k \varphi\}_{k=-\infty}^\infty$
with a bounded operator $T,$ but it is an open problem whether it also
has a representation with the indexing in \eqref{70404a}. For Gabor systems it is
not known whether the chosen indexing play a role.

The question of the existence of a representation on the form \eqref{70404a} is
interesting for general frames as well. We will also formulate an open problem
for general frames in an arbitrary separable Hilbert space; a negative answer to
that question would also lead to a negative answer to the above specific questions.

In Section \ref{70331p} we will provide the necessary background on frames and
operator representations. The presentation is kept as short as possible, and
only contains the information that is necessary in order to understand the open
problems; these  are outlined in Section \ref{70403a}.
Section \ref{70331s} contains further motivation, analysis of special cases, and results about certain orderings
of the frame elements that must be avoided if we want to obtain a representation
of the  form \eqref{70404a}.

\section{Technical background} \label{70331p}

\subsection{Frame theory} Let $\h$ denote a separable Hilbert space. A sequence $\ftki$
in  $\h$ is a {\it frame}
for  $\h$ if there exist constants $A,B>0$ such that
$ A \, ||f||^2 \le \sumk | \la f, f_k\ra|^2 \le B \, ||f||^2, \, \forall f\in \h;$ it is
a {\it frame sequence} if the stated inequalities hold for all $f\in \span \ftki.$
The sequence $\ftki$ is a {\it Bessel sequence} if at least the upper frame condition holds. Also, $\ftki$ is called a {\it Riesz sequence} if
there exist constants $A,B>0$ such that
$A \sum |c_k|^2 \le \nl \sum c_k f_k \nr^2 \le B\sum |c_k|^2$
for all finite scalar sequences $\{c_k\}_{k\in I}.$ A {\it Riesz basis} is a Riesz
sequence $\ftki$ for which $\span \ftki = \h.$

If $\ftki$ is a Bessel sequence, the {\it synthesis operator} is defined by
\bee \label{60811a} U: \ell^2(I) \to \h, \, U\ctki := \sumk c_k f_k;\ene
it is well known that $U$ is well-defined and bounded. A central role will be played by the kernel
of the operator $U,$ i.e., the subset of $\lti$ given by
\bee \label{60811f} N_U=\left\{\ctki \in \ell^2(I) ~\bigg|~\sumk c_kf_k=0\right\}.\ene
The {\it excess} of a frame  is the number of elements that can be removed yet leaving a frame.
It is well-known that the excess equals $\mbox{dim} (N_U);$ see  \cite{Balan}.

Given a Bessel sequence $\ftki,$ the {\it frame operator} $S: \h \to \h$ is defined by $S:=UU^*.$
For a frame $\ftki,$ the frame operator is invertible and
$f= \sum_{k\in I} \la f, \si f_k\ra f_k, \, \forall f\in \h.$
The sequence $\{\si f_k\}_{k\in I}$ is also a frame; it is called the {\it canonical dual frame.}

Note that every orthonormal basis (or more generally, every Riesz basis) $\ftki$
is a frame.  A frame that is not a basis is said to be overcomplete. We refer to
\cite{CB} and \cite{Heil} for more information about frames and Riesz bases.

\subsection{Operator representations of frames}

Dynamical sampling typically concern frame properties of sequences in a Hilbert space
$\h$ of the form
$\{T^k \varphi\}_{k=0}^\infty,$ where $\varphi \in \h$ and
$T: \h \to \h$ is a linear operator. We will also consider sequences indexed by $\mz;$
thus, in the sequel the index set $I$ denotes  either $\mz$ or $\mn_0= \{0\} \cup \mn.$

As in the
papers \cite{olemmaarzieh-1,olemmaarzieh-2} we will look at the topic from the opposite
side.  Indeed, we will consider a given  frame $\ftki$ in a Hilbert space $\h$
and ask for the existence of a representation on the
form
\bee \label{61209a} \ftki = \{T^k \varphi\}_{k\in I}, \ene
where $T: \Span \ftki  \to \h$ is a bounded linear operator. Note that we will consider
$\ftki$ as a sequence, i.e., as an ordered set.

The following result collects
findings from the papers \cite{olemmaarzieh-1, olemmaarzieh-2}, describing when
a representation of the form \eqref{61209a} is possible, and when the operator $T$
 can be chosen to be bounded. We will need the right-shift operator on $\lti,$ defined by
\bes {\cal T}: \lti \to \lti, {\cal T}\ctki= \{c_{k-1}\}_{k\in I};\ens
for the case $I= \mn_0$ we define $c_{-1}:=0.$

\bpr \label{511p} Consider a frame $\ftki$ for an
infinite-dimensional Hilbert space $\h,$ with $I$ denoting either
$\mn_0$ or $\mz$. Then the following hold:
\bei \item[(i)] There exists a linear operator $T: \mbox{span} \ftki \to \h$
such that \eqref{61209a} holds if and only if $\ftki$ is linearly independent.
\item[(ii)] Assume that $\ftki$ is linearly independent. Then
the operator $T$ in \eqref{61209a} is bounded if and only the kernel $N_U$
of the synthesis operator is invariant under right-shifts; in particular $T$ is bounded if $\ftki$ is a Riesz basis.
\item[(iii)] Assume that $\ftki$ is linearly independent
and overcomplete. If the operator $T$ in \eqref{61209a} is bounded, then $\ftki$ has infinite excess.
\eni \epr

Let us mention an example of an overcomplete frame that indeed has a representation
on the form \eqref{70404a}; the example is due to Aldroubi et al. \cite{A2}.

\bex \label{70422a} Consider  the
matrix $T= [a_{ij}]_{i,j\in \mn}$ given by $a_{jj}= 1-2^{-j}, \, a_{ij}= 0, \, i\neq j,$
and the sequence $g:= \{ \sqrt{1-(1-2^{-j})^2}\}_{j\in \mn}.$ Then the system $\{T^kg\}_{k=0}^\infty$ is a frame for $\ltn$
but not a basis. 
\ep \enx

\subsection{Shift-invariant systems and Gabor frames}

Shift-invariant systems and Gabor systems are defined in terms of certain classes
of operators on $\ltr.$
For $a\in \mr,$ define the {\it translation operator}
$T_a$ acting on $\ltr$ by $T_af(x):= f(x-a)$ and the {\it modulation operator}
$E_a$ by $E_af(x):=e^{2\pi i ax} f(x).$  Both operators are unitary.  Furthermore,
defining the Fourier transform of $f\in L^1(\mr)$ by $\widehat{f}(\ga)={\cal F}f(\ga)=
\inr f(x)e^{-2\pi i \ga x}dx$ and extending it in the standard way to a unitary operator
on $\ltr,$ we have ${\cal F}T_a=E_{-a}{\cal F}.$

Given a function $\varphi\in \ltr$
and some $b>0$,  the associated {\it shift-invariant system} is given by
$\{T_{kb}\varphi\}_{k\in\mz}.$ The following result collects the necessary information about such systems. Consider the function

\bee \label{70331a} \Phi(\ga):= \sum_{k\in \mz}|\widehat{\varphi}(\frac{\ga + k}{b})|^2, \, \ga \in \mr.\ene

\bpr \label{70331q} Let $\varphi \in \ltr \setminus \{0\},$ and $b>0$ be given. Then the following hold:
\bei
\item[(i)] $\{T_{kb}\varphi\}_{k\in\mz} $ is linearly independent.
\item[(ii)] $\{T_{kb}\varphi\}_{k\in\mz} $ is a Riesz basis if and only if
there exist $A,B>0$ such that $A\le\Phi(\ga)\le B,$ a.e. $\ga\in [0,1].$
\item[(iii)] $ \{T_{kb}\varphi\}_{k\in\mz} $ is a frame sequence if and only if
there exist $A,B>0$ such that $A\le\Phi(\ga)\le B,$ a.e. $\ga\in [0,1]\setminus \{ \ga\in [0,1]~ \big| ~  \Phi(\ga) =0    \}.$
\item[(iv)] If $ \{T_{kb}\varphi\}_{k\in\mz} $ is an overcomplete frame sequence, it has
infinite excess.
\item[(v)]
$\{T_{kb}\varphi\}_{k\in\mz}= \{(T_b)^k \varphi\}_{k\in\mz}, $ i.e.,  the system $\{T_{kb}\varphi\}_{k\in\mz} $
has the form of an iterated system indexed by $\mz.$
\eni
\epr
The result in (i) is well-known, and (ii) \& (iii) are proved in \cite{BL1}; (iv)
is proved in \cite{Balan,olemmaarzieh-2}, and (v) is clear.

\bex \label{70404b} Letting
\bee \label{203ag}
\mbox{sinc}(x):=\left\{
\begin{array}{lll}
{\large \frac{\sin(\pi x)}{\pi x}} & \mbox{if  $ x \neq 0$;} \\
1 & \mbox{if $ x =0$,} \\
\end{array}
\right. \ene 
Shannon's sampling theorem states that the shift-invariant system
$\{T_k \mbox{sinc}\}_{k\in \mz}$ \ form an
orthonormal basis for the {\it Paley-Wiener space}
\bee \label{507a} PW:= \left\{f\in \ltr \
\big| \ \supp \, \widehat{f} \subseteq \big[-\frac12,\frac12
\big]\right\}.\ene
It follows from this that the system $\{T_{k/2} \mbox{sinc}\}_{k\in \mz}$ is
an overcomplete frame for the Paley-Wiener space; indeed,
$\{T_{k/2} \mbox{sinc}\}_{k\in \mz}=
\{T_{k} \mbox{sinc}\}_{k\in \mz} \bigcup \{T_{1/2} T_{k} \mbox{sinc}\}_{k\in \mz},$
which shows that the system is a union of two orthonormal bases for the
Paley-Wiener space.
\ep \enx

A collection of functions in $\ltr$ of the form $\mts$ for some $a,b>0$ and some
$g\in \ltr$ is called a {\it Gabor system}. Gabor systems play an important role
in time-frequency analysis; we will just state the properties that are necessary
for the flow of the current paper, and refer to \cite{G2,FS1,FS2,CB} for
much more information.

\bpr \label{61213a} Let $g\in \ltr \setminus \{0\}.$ Then the following hold:

\bei
\item[(i)] $\mts$ is linearly independent.
\item[(ii)] If $\mts$ is a frame for $\ltr,$ then $ab\le 1.$
\item[(iii)]  If $\mts$ is a frame for $\ltr,$ then $\mts$ is a Riesz basis if and only if $ab= 1.$
\item[(iv)]  If $\mts$ is an overcomplete frame for $\ltr$, then $\mts$ has infinite excess.
\eni
\epr
The result in (i) was proved in \cite{Linn}  (hereby confirming a conjecture stated in
\cite{HRT});  (ii) \& (iii) are classical results \cite{G2,CB}, and (iv) is proved in \cite{Balan}.

Note that it immediately follows from Proposition \ref{511p} and Proposition \ref{61213a} that any
Gabor frame $\mts$  can be represented on the form
$\Tni$ for some linear operator $T:\ltr \to \ltr$ and some $\varphi \in \ltr.$
Whether there exists a {\it bounded} operator $T$  with this property 
is a much harder question, and indeed the topic for the next section. Let us apply Proposition \ref{511p} and show an example of a concrete
Gabor frame and a concrete ordering that leads to an unbounded operator.

\bex \label{70404u} Consider the Gabor frame  $\{E_{m/3}T_{n}\chi_{[0,1]}\}_{m,n\in\mz},$ which is
the union of the three orthonormal bases
$\{E_{k/3}E_{m}T_{n}\chi_{[0,1]}\}_{m,n\in\mz}, \, k=0,1,2.$ Re-order the frame as $\ftz$ in such a way that the elements $\{f_{2k+1}\}_{k\in\mz}$ corresponds to the orthonormal basis $\{E_{m}T_{n}\chi_{[0,1]}\}_{m,n\in\mz}.$ By construction, the elements $\{f_{2k}\}_{k\in\mz}$ now forms an overcomplete frame. By Proposition \ref{511p}, there is an operator $ T:\Span\ftz\rightarrow \Span\ftz $ such  that $ \ftz = \Tnz $. Since the subsequence $ \{f_{2k}\}_{k\in\mz} $ is an overcomplete frame, there is a non-zero sequence $\{ c_{2k} \}_{k\in\mz}\in \ltz $ such that $ \sum_{k\in\mz}c_{2k}f_{2k}=0 $. Defining $c_k=0 $ for $ k\in 2\mz + 1 $, we have $\sum_{k\in\mz}c_{k}f_{k}=\sum_{k\in\mz}c_{2k}f_{2k}=0$.
On the other hand, since   $ \{f_{2k+1}\}_{k\in\mz} $   is a Riesz basis and $ \ctz $ is non-zero, $ \sum_{k\in\mz}c_{k}f_{k+1}=\sum_{k\in\mz}c_{2k}f_{2k+1} \neq 0 $.  This shows that $N_U$ is not invariant under right-shifts; thus, $T$ is unbounded by Proposition \ref{511p}.
\ep \enx

\section{The open problems} \label{70403a}
In this section we continue our convention of letting the index set $I$
denote either $\mn_0$ or $\mz.$
The results in Section \ref{70331p} show that there is a perfect match between the
stated necessary conditions for a frame to be represented in terms of a bounded operator
as in \eqref{61209a},
and the properties of Gabor frames and shift-invariant systems. Indeed, by Proposition \ref{511p} a frame for an infinite-dimensional Hilbert space has the form
\eqref{61209a} for a (not necessarily bounded) operator $T$ if it is
linearly independent; by Proposition \ref{70331q} and Proposition \ref{61213a} this condition
is satisfied for
the Gabor frames and the shift-invariant frames. Also, Proposition \ref{511p} shows that if
a frame $\ftki$ for $\h$ is linearly independent and overcomplete, then the operator
$T$ can only be bounded if $\ftki$ has infinite
excess; again, by Proposition \ref{70331q} and Proposition \ref{61213a} these conditions  match
the Gabor frames and the shift-invariant frames. For the case of a shift-invariant frame
we have already seen in Proposition \ref{70331q}  that it has the form $\{(T_b)^k \varphi\}_{k\in \mz},$ i.e., as an iterated system indexed by $\mz.$
We will therefore first state the open problems for Gabor frames, and return to
shift-invariant systems afterwards.

On short form, the key problem is as follows:

\vn {\bf Problem 1:} Does there exist overcomplete Gabor frames $\mts$
and an appropriate ordering on the form $\ftki,$ such that the
operator $T$ in \eqref{61209a} is bounded?

\vspace{.5cm}  From the formulation of Problem 1 it is already clear that the
question contains several subproblems. Assuming that Problem 1 has a positive answer
for a certain overcomplete Gabor frame $\mts,$
we can formulate them as follows:
\bei
 \item[(P1-a)] Is a representation of the form \eqref{61209a} with $T$ bounded available
for both indexings $I=\mn_0$ and $I=\mz$?
\item[(P1-b)] Can we characterize the orderings $\ftki$ of $\mts$ for which a
representation as in \eqref{61209a} is available with a bounded operator $T$?
\item[(P1-c)] Can we characterize the  Gabor frames $\mts$ for which
an appropriate ordering on the form $\ftki$ leads to a bounded
operator $T$ in \eqref{61209a}, i.e., what are the conditions
on $a,b>0$ and $g\in \ltr$ for this to be possible?
\item[(P1-d)] Is there  a universal
ordering on the form $\ftki$ that applies for all the Gabor frames in (P1-c)?
\eni
As far as we know, no information is available in the literature concerning
question (P1-a).  Concerning (P1-b), recall that an ordering that lead to
a negative conclusion was considered in Example \ref{70404u}.

Let us now turn our attention to shift-invariant systems $\{T_{kb}\varphi\}_{k\in \mz}.$
Using Proposition \ref{70331q} (v) we immediately see that
such a system always has the representation $\{(T_b)^k\varphi\}_{k\in \mz},$ with
the translation operator $T_b$ being bounded. Thus Problem 1
trivially has a positive answer  for the index set $I=\mz.$ However, it
is unclear whether the same is true for the index set $I=\mn_0.$  Let us
formulate the key question:

\vn {\bf Problem 2:} Is it possible to order an overcomplete frame of translates
$\{T_{kb}\varphi\}_{k\in \mz}$ on the form $\{f_k\}_{k=0}^\infty = \{T^k \varphi\}_{k=0}^\infty$ for a bounded operator $T$?

\vn For the question of an operator representation indexed by
$\mn_0,$ there is an even more fundamental problem underlying Problem 1 \& 2. 
The problem can be phrased in the setting of a general Hilbert space.
Indeed,
the known cases of overcomplete frames
having a representation $\{f_k\}_{k=0}^\infty =\{T^k f_0\}_{k=0}^\infty$ in terms of a bounded operator have the property that $f_k \to 0$ as $k\to \infty;$ see, e.g., Example \ref{70422a}.

\vn {\bf Problem 3:} Let $\h$ denote a separable Hilbert space. 
Do there exist overcomplete frames for $\h$ that are norm-bounded below
and have a representation $\{f_k\}_{k=0}^\infty = \{T^k f_0\}_{k=0}^\infty$ for a bounded operator $T:\h \to \h$?

\vn Note that since Gabor frames and
shift-invariant systems consist of vectors with
equal norm, a negative answer to Problem 3 would automatically lead to a negative
answer to Problem 1  (for the index set $I=\mn_0)$ and Problem 2.

\section{Frames in Hilbert spaces}  \label{70331s}

The purpose of this section is to shed light on the questions in Section \ref{70403a}
from a general Hilbert space angle. We first show that 
an overcomplete frame with a representation on the form \eqref{70404a} must have
a very particular property:

\bl \label{70428b} Assume that $\{f_k\}_{k=0}^\infty$ is an overcomplete frame and that 
$\{f_k\}_{k=0}^\infty= \{T^kf_0\}_{k=0}^\infty$ for some bounded linear
operator $T: \h \to \h.$ Then there exists an $N\in \mn_0$ such that
$\{f_k\}_{k=0}^N \cup \{f_k\}_{k=M}^\infty$ is a frame for $\h$ for all $M>N.$ \el

\bp Choose
some coefficients $\{c_k\}_{k=0}^\infty\in \ell^2(\mn_0)$ such that
$\sum_{k=0}^\infty c_kf_k=0.$ Letting $N:= \min\{k\in \mn_0 \, | \, c_k \neq 0\},$
we have that 
\bee \label{70428a} -c_Nf_N= \sum_{k=N+1}^\infty c_kf_k,\ene so
$f_N\in \span \{f_k\}_{k=N+1}^\infty.$ Thus $\{f_k\}_{k=N+1}^\infty$
is a frame for $\span \{f_k\}_{k=N}^\infty.$ Applying 
the operator $T$ on \eqref{70428a} shows that $f_{N+1}\in \span \{f_k\}_{k=N+2}^\infty.$ By
iterated application of the operator $T$ this proves that for any $M>N,$ the family $\{f_k\}_{k=M}^\infty$ is
a frame for $\span \{f_k\}_{k=N}^\infty,$ which leads to the desired result. \ep

Intuitively, Lemma \ref{70428b} says that a frame having the desired type of representation
must have ``infinite excess in all directions."

We will now consider frames with a special structure
and examine the existence of a representation on the form \eqref{70404a}.
Remember from Proposition \ref{511p} that any orthonormal basis has a representation
on the form \eqref{70404a}; on the other hand, for a sequence consisting
 of a union of a basis and finitely many elements, a representation on this form is
not possible with a  bounded operator $T$. In order to analyse the existence of
a representation \eqref{70404a} it is therefore natural to consider the union of
two orthonormal bases; the indexing of these orthonormal bases is not relevant, and
we will write them as $\etk, $ respectively, $\eptk.$ Note that this setup
actually covers certain special Gabor frames, and hereby directly relate to
the questions in Section \ref{70403a}; for example, the Gabor frame $\{E_{m/2}T_n \chi_{[0,1]}\}_{m,n\in \mz}$ is indeed the union of the two orthonormal bases
$\{E_{m}T_n \chi_{[0,1]}\}_{m,n\in \mz}$ and $\{E_{1/2}E_mT_n \chi_{[0,1]}\}_{m,n\in \mz}.$

Thus, let us now consider a frame $\ftk$ in a Hilbert space $\h,$ formed as the union of the elements in two
orthonormal bases $\etk$ and $\eptk.$ In order to describe the
elements in $\ftk$ let us introduce the index sets
\bes I_1 & :=&  \{k\in \mn \, \big| \, f_k \in \{e_j\}_{j=1}^\infty  \mbox{ and }
f_{k+1} \in \{\varepsilon_j\}_{j=1}^\infty\}; \\
I_2 & :=&  \{k\in \mn \, \big| \, f_k \in \{e_j\}_{j=1}^\infty  \mbox{ and }
f_{k+1} \in \{e_j\}_{j=1}^\infty\}; \\
I_3 & :=&  \{k\in \mn \, \big| \, f_k \in \{\varepsilon_j\}_{j=1}^\infty  \mbox{ and }
f_{k+1} \in \{e_j\}_{j=1}^\infty\}; \\
I_4 & :=&  \{k\in \mn \, \big| \, f_k \in \{\varepsilon_j\}_{j=1}^\infty  \mbox{ and }
f_{k+1} \in \{\varepsilon_j\}_{j=1}^\infty\}.
\ens
Note that the sets $I_k, k=1, \dots,4$ form a disjoint covering of the index set $\mn$
of the frame $\ftk,$ and that the sets $I_1$ and $I_3$ always are nonempty and infinite sets.
By symmetry we can assume that $f_1=e_1.$ The following result shows that
{\it if} it is possible for the frame $\ftk$ to have a representation
as in \eqref{70404a} with $I=\mn_0,$ then necessarily $I_2 \neq \emptyset$ and $I_4 \neq \emptyset.$

\bl \label{70422d} In the setup described above, assume that either $I_2= \emptyset$ or
$I_4= \emptyset.$ Then $T$ is unbounded.
\el

\bp Assume that $\ftk$ has a representation on the
form \eqref{70404a}, with $T$ being bounded. Then $T$ is surjective.  Indeed, since
$f_1=e_1,$ then the set $\{Tf_k\}_{k=1}^\infty$ contains all the vectors $\eptk.$
This implies
that $T$ is surjective.

Now assume that $I_2 \neq \emptyset$ and $I_4= \emptyset.$ Then
$\mn = I_1 \cup I_2 \cup I_3,$ and $\{T \varepsilon_k\}_{k=1}^\infty \subseteq \{e_k\}_{k=2}^\infty.$ Since $\eptk$ is an orthonormal basis for $\h,$ this implies that $\cra(T) \subseteq
\span\{e_k\}_{k=2}^\infty \neq \h.$ This contradicts the surjectivity of
$T.$ 
Thus the operator $T$ can not be bounded in this case. The case where 
$I_2= \emptyset$ and $I_4 \neq \emptyset$ is similar. In this case there
exist some $\ell, \ell^\prime \in \mn$ such that $T\varepsilon_\ell = \varepsilon_\ell^\prime,$ and thus $\{T e_k\}_{k=1}^\infty \subseteq \{\varepsilon_k\}_{k\neq \ell^\prime}.$ Again this contradicts the surjectivity of $T.$ Thus, again
$T$ can not be bounded. \ep

In particular, Lemma \ref{70422d} shows that the operator $T$ can not
be bounded if the elements in the frame $\etk \cup \eptk$ are ordered as
$\{e_1, \varepsilon_1, e_2, \varepsilon_2, \dots\}.$ This shows
a significant difference between the availability of a representation
as in \eqref{70404a} for $I= \mn_0$ and $I=\mz:$ indeed, Example \ref{70404b}
yields a concrete case where a union of two orthonormal bases, ordered as
$\{\cdots,\varepsilon_0, e_0, \varepsilon_1, e_1,\varepsilon_2,e_2,\cdots\},$ actually
has a representation as in \eqref{70404a} with $I=\mz.$
This difference between the indexing in terms of $I=\mn_0$ and $I=\mz$ is
precisely our motivation for Problem (P1-a) and Problem 2 from the
previous section.

\begin{tabbing}
text-text-text-text-text-text-text-text-text-text \= text \kill \\
Ole Christensen \> Marzieh Hasannasab \\
Technical University of Denmark \> Technical University of Denmark  \\
DTU Compute \> DTU Compute \\
Building 303, 2800 Lyngby \> Building 303, 2800 Lyngby \\
Denmark \> Denmark \\
Email: ochr@dtu.dk \> mhas@dtu.dk
\end{tabbing}


\begin{thebibliography}{99}

\bibitem{A1} Aldroubi, A.,  Cabrelli, C.,  Molter, U., and Tang, S.:
{\it Dynamical sampling.} Appl. Harm. Anal. Appl. {\bf 42} no. 3 (2017),  378-–401.


\bibitem{A2} Aldroubi, A., Cabrelli, C., \c{C}akmak, A. F., Molter, U., and  Petrosyan, A.: {\it Iterative actions of normal operators.} J. Funct. Anal. {\bf 272} no. 3 (2017),
    1121--1146.

\bibitem{A3} Aldroubi, A., and Petrosyan, A.: {\it Dynamical sampling and systems
from iterative actions of operators.} Preprint, 2016.

\bibitem{Balan} Balan, R., Casazza, P., Heil, C. and Landau, Z.: {\it Deficits
and excesses of frames.} Adv. Comp. Math. {\bf 18} (2002), 93--116.

 \bibitem{BL1} Benedetto, J. and Li, S.:
{\it The theory of multiresolution analysis frames and
applications to filter banks.} Appl. Comp. Harm. Anal., {\bf 5}
(1998), 389--427.



\bibitem{CB} Christensen, O.:
{\it An introduction to frames and Riesz bases.} Second expanded
edition. Birkh\"auser (2016)





\bibitem{olemmaarzieh-1}
		Christensen, O., and Hasannasab, M.:
		\newblock {\em Frame properties of systems arising via iterative actions of operators}
		\newblock Preprint, 2016.

\bibitem{olemmaarzieh-2}
		Christensen, O., and Hasannasab, M.:
		\newblock {\em Operator representations of frames: boundedness, duality, and stability}
		\newblock Accepted for publication in Integral Equations an Operator Theory.

\bibitem{Da1}  Daubechies, I.: {\it The wavelet transformation,
time-frequency localization and signal analysis.} IEEE Trans.
Inform. Theory {\bf 36} (1990), 961--1005.



\bibitem{FS1} Feichtinger, H. G. and Strohmer, T. (eds.): {\it Gabor Analysis
and Algorithms: Theory and Applications.} Birkh\"auser, Boston,
1998.

\bibitem{FS2} Feichtinger, H. G. and Strohmer, T. (eds.): {\it
Advances in Gabor Analysis.} Birkh\"auser, Boston, 2002.

\bibitem{G2} Gr\"{o}chenig, K.: {\it Foundations of
time-frequency analysis.} Birkh\"{a}user, Boston, 2000.

\bibitem{Heil} Heil, C.: {\it A basis theory primer.} Expanded edition. Applied and Numerical Harmonic Analysis. Birkhäuser/Springer, New York, 2011.



\bibitem{HRT} Heil, C., Ramanathan, J. and Topiwala, P.: {\it Linear
independence of time-frequency translates.} Proc. Amer. Math. Soc.
{\bf 124} (1996),  2787--2795.

\bibitem{Klauder} Klauder, J. and Skagerstam, B.: {\it Coherent states.  Applications in physics and mathematical physics.} World Scientific, 1985.



\bibitem{Linn} Linnell, P.: {\it Von Neumann algebras and linear independence
of translates.} Proc. Amer. Math. Soc. {\bf 127} no. 11 (1999),
3269--3277.

\bibitem{Perelomov} Perelomov, A.: {\it Generalized Coherent States and Their Applications.}
Springer 2012.

\bibitem{FP} Philipp, F.: {\it Bessel sequences from iterated operator actions.} Preprint, 2016.

\end{thebibliography}
\end{document}